\documentclass[11pt]{article}
\usepackage{amsmath}
\usepackage{amsthm}
\usepackage{amssymb}
\usepackage{graphics}
\usepackage{enumerate}   
\newcommand{\be}{\begin{equation}}
\newcommand{\ee}{\end{equation}}
\newcommand{\ben}{\begin{eqnarray*}}
\newcommand{\een}{\end{eqnarray*}}
\newtheorem{examp}{\sc Example}
\newtheorem{remk}{\sc Remark}
\newtheorem{corol}{\sc Corollary}
\newtheorem{lemma}{\sc Lemma}
\newtheorem{theorem}{\sc Theorem}
\newtheorem{defn}{\sc Definition}
\newcommand{\bt}{\begin{theorem}}
\newcommand{\et}{\end{theorem}}
\newcommand{\bl}{\begin{lemma}}
\newcommand{\el}{\end{lemma}}
\newcommand{\bed}{\begin{defn}}
\newcommand{\eed}{\end{defn}}
\newcommand{\brem}{\begin{remk}}
\newcommand{\erem}{\end{remk}}
\newcommand{\bex}{\begin{examp}}
\newcommand{\eex}{\end{examp}}
\newcommand{\bcl}{\begin{corol}}
\newcommand{\ecl}{\end{corol}}

\newcommand{\NI}{\noindent}
\newcommand{\lmd}{\lambda}
\newcommand{\al}{\alpha}

\newcommand{\fal}{\forall}
\newcommand{\raro}{\rightarrow}

\newcommand{\dsp}{\displaystyle}
\newcommand{\vsp}{\vskip 0.5em}

\newcommand{\cs}{{\cal S}}
\newcommand{\ct}{{\cal T}}

\theoremstyle{definition}
\theoremstyle{remark}

\numberwithin{equation}{section}
\numberwithin{theorem}{section}
\numberwithin{lemma}{section}

\begin{document}
\title{Some Aspects on Solving Transportation Problem}
\author{A. K. Das$^{a,1}$, Deepmala$^{a,2}$, and R. Jana$^{b,3}$\\
	\emph{\small $^{a}$Indian Statistical Institute, 203 B. T.
		Road, Kolkata, 700 108, India.}\\
	\emph{\small $^{b}$Jadavpur University, Kolkata , 700 032, India.}\\
	\emph{\small Email: $^{1}$akdas@isical.ac.in, $^{2}$dmrai23@gmail.com, $^{3}$rwitamjanaju@gmail.com.}}
\date{}
\maketitle

\begin{abstract}
	\noindent In this paper, we consider a class of transportation problems which arises in sample surveys and other areas of statistics. The associated cost matrices of these transportation problems are of special structure. We observe that the optimality of North West corner solution holds for the general problem where cost component is replaced by a convex function. We revisit assignment problem and present a weighted version of K$\ddot{\mbox{o}}$nig-Egerv$\acute{\mbox{a}}$ry theorem and Hungarian method. The weighted Hungarian method proposed in the paper can be used for solving transportation problem.\\
	
	\NI{\bf Keywords: } Transportation problem; North West corner solution; weighted K$\ddot{\mbox{o}}$nig-Egerv$\acute{\mbox{a}}$ry theorem;  assignment problem; weighted Hungarian  method; sample survey
	
\end{abstract}

\footnotetext[2]{Corresponding author}


\section{Introduction}
\noindent A transportation model is a bipartite graph $G=(A\cup B, E)$ where $A$ $=\{O_{1},\cdots,O_{m}\}$ is the set of source vertices, $B$ $=\{D_{1},\cdots,D_{n}\}$ is the set of destination vertices and $E$ is the set of edges from $A$ to $B.$ Each edge $(i,j)\in E$ has an associated cost $c_{ij}.$ The problem is to find out a flow of least costs that ships from supply sources $O_{i},$ $i=1,\cdots, m$ to consumer destinations $D_{j},$ $j=1,\cdots, n.$ Suppose $a_{i}$ is the supply at the $i^{th}$ source $O_{i}$ and $b_{j}$ is the demand at the $j^{th}$ destination $D_{j}.$ In a balanced transportation problem, we assume that $\dsp {\sum_{i=1}^{m}a_{i}=
	\sum_{j=1}^{n} b_{j}}.$ Let $x_{ij}$ be the quantity to be shipped from origin
$O_{i}$ to destination $D_{j}$ with cost  $c_{ij}.$ The transportation problem can be formulated as a linear programming problem to determine a shipping schedule
that minimizes the total cost of shipment which is given below.\\

$
\begin{array}{lllll}
& \quad \quad \;  &     & \mbox{minimize} &\dsp{\sum_{i\in A} \sum_{j\in B}}c_{ij}x_{ij}\\
& & & \mbox{subject to} & \\
& & & & \dsp{\sum_{j\in B}}x_{ij}=a_{i},\;\fal\;i\in A\\
& & & & \dsp{\sum_{i\in A}}x_{ij}=b_{j},\;\fal\;j\in B\\
& & & & \dsp{x_{ij}\geq 0,\;i\in A,\;j\in B.}
\end{array}
$

\vsp
\noindent An assignment problem is a special case of a balanced transportation problem where $m=n,$  $a_{i}=1,\;\fal\; i\in A$ and  $b_{j}=1,\;\fal\;j\in B.$  Various generalizations of transportation problem  have been appeared in the literature. For details, see Goossens and Spieksma \cite{1},  Kasana and Kumar \cite{kasana},  Liu and Zhang \cite{11} and the references cited therein.

\noindent We show that a class of transportation problem arises in statistics. We discuss the various structures and solution method of this class. We also consider the Hungarian method for the assignment problem and extend it for solving a transportation problem. In section 2, we consider the class of transportation problem and its application in statistics.  We present an elegant proof of the result that the North West Corner rule provides an optimal solution to the transportation problem under some conditions. In section 3, we consider assignment problem and present a weighted version of K$\ddot{\mbox{o}}$nig-Egerv$\acute{\mbox{a}}$ry theorem and Hungarian method. We show  minimum cut - maximum flow theorem of Ford and Fulkerson in a different way so that in a bipartite graph, finding a minimum weight vertex cover is equivalent with finding a minimum cut if the capacity of an edge is given by the minimum weight of its end nodes. The weighted version of Hungarian method proposed in the paper for solving transportation problem is same as primal-dual method for solving minimum cost network flow problems. We establish a connection between transportation problem and assignment problem and propose a weighted version of Hungarian method to solve transportation problem. Section 4 presents the conclusions of the paper.

\section{Transportation Problem and its Applications in Statistics}
Transportation problem arises in various applications of Sample Surveys and Statistics. The structure of the cost matrices associated with these transportation problems are of special structure. Now, we raise the following question. What are the structure of the cost matrix for which North West corner solution produces an optimal solution? We consider some of the structures of the cost matrix which arise in some of the applications in the literature.  Hoffman \cite{15} studied transportation problem in the context of North West Corner Rule.  Burkard et al. \cite{14} mentioned Monge properties in connection with   transportation problem. Szwarc \cite{17}  developed direct methods for solving transportation problems with cost coefficients of the form $c_{ij}= x_i + x_j$, having applications in shop loading and aggregate scheduling.

Evans \cite{6} studied factored transportation  problem in which cost coefficients are factorable, i.e., $c_{ij}= x_i x_j$.  It is shown that   the rows and columns can easily be ordered so that the North West corner rule provides an optimal solution of the transportation problem.  We state some of the results of Evans \cite{6} which are needed in the sequel.

\bl
The North West corner rule produces an optimal solution of the balanced transportation problem whenever $c_{ij} + c_{rs} \leq c_{rj} + c_{is}$ for all $i,j, r, s$ such that $i < r$ and $j < s.$
\el

\bt
Let $x_{i},\; i = 1,2,$ $\cdots, m$ and $y_{j},\;j = 1,2,\cdots,n$ be nonnegative real numbers such that $x_{1}\geq x_{2}>\cdots \geq  x_m$ and $y_1\leq y_2<\cdots \leq y_n.$ Then the North West corner rule provides an optimal solution to the transportation problem where $c_{ij} = x_i y_j.$
\et

\bt
Let $x_{i},\; i = 1,2,\cdots, m$ and $y_{j},$ $j = 1,2,\cdots,n$ be  real numbers.   Then the North West corner rule provides an optimal solution to the transportation problem where $c_{ij} = x_i + y_j.$
\et

Szwarc \cite{17} showed that if $c_{ij} = x_i +y_j$ then any feasible solution is optimal.
Raj \cite{16} studied the problem of integration of surveys, i.e., the problem of designing a sampling program for two or more surveys which maximizes the overlap between observed samples as a transportation problem and shows that the solution is just the North West corner solution of the transportation problem and this is optimal when $c_{ij}=|i-j|.$   For the connection between integration of surveys and the transportation problem  see Matei and Till$\acute{\mbox{e}}$ \cite{12}, Aragon and Pathak \cite{3},  Causey et al. \cite{2} and Raj \cite{16}. In this context, Burkard et al. \cite{14} studied several perpectives of Monge properties in optimization. Mitra and Mohan \cite{13} observed that the North West corner solution is optimal  for the following problem.

Suppose $X$ and $Y$ are two discrete random variables which assume values $x_{1}\leq x_{2}\leq \cdots \leq x_{m}$ and $y_{1}\leq y_{2}\leq \cdots\leq y_{n}$ respectively. Let $p_{i\cdot}= \mbox{Prob} (X=x_{i})$ and $p_{\cdot j}= \mbox{Prob} (Y=y_{j}).$ The problem is to find out the joint probabilities $p_{ij}= Prob(X=x_{i}, Y=y_{j})$ so that $\mbox{Cov}(X,Y)$ is maximum. This problem can be formulated as a transportation problem as follows.

Given values of random variables $X,$  $Y,$  $p_{i\cdot},$ $p_{\cdot j}, $ the problem is to find $p_{ij},$ $1\leq i\leq m, $ $1\leq j\leq n$ which

$
\begin{array}{lllll}
& \quad \quad \;  &     & \mbox{minimize} &\dsp{ \sum_{i}\sum_{j} (x_{i}-y_{j})^{2}\,p_{ij}}\\
& & & \mbox{subject to} & \\
& & & & \dsp{ \sum_{j}p_{ij}= p_{i\cdot},\;1\leq i\leq m}\\
& & & & \dsp{ \sum_{i}p_{ij}= p_{\cdot j},\;1\leq j\leq n}\\
& & & & \dsp{ p_{ij}\geq 0,\;1\leq i\leq m,\;1\leq j\leq n}
\end{array}
$

\vsp

In this article, we consider a more general problem and show that the optimality of North West corner solution holds.

Suppose $f:\;R\raro R$ is a convex function. Let $\dsp{ c_{ij}= f(x_{i}-y_{j})}.$
Given values of random varibles $X$,  $Y$,  $p_{i\cdot}, p_{\cdot j}, $ the problem is to find $p_{ij}$, $1\leq i\leq m,$ $1\leq j\leq n$ for the Problem P which is stated as follows.\\

$
\begin{array}{lllll}
&\bf Problem\ P: \quad \quad \;  &     & \mbox{minimize} &\dsp{ \sum_{i}\sum_{j} f(x_{i}-y_{j})\,p_{ij}}\\
& & & \mbox{subject to} & \\
& & & & \dsp{ \sum_{j}p_{ij}= p_{i\cdot},\;1\leq i\leq m}\\
& & & & \dsp{ \sum_{i}p_{ij}= p_{\cdot j},\;1\leq j\leq n}\\
& & & & \dsp{ p_{ij}\geq 0,\;1\leq i\leq m,\;1\leq j\leq n}
\end{array}
$

\vsp
Let us define a set $B=\{(i,j)\;|\; p_{ij}\;\mbox{ is a basic variable}\}.$ Note that if $B$ is the basis set and if the basic solution to P corresponding to the set $B$ is also feasible, it is optimal if and only if there exist $\al_i$, $1\leq i\leq m$, $\beta_j$, $1\leq j\leq n$ such that
\be
\al_i +\beta_j = c_{ij}\;\mbox{ if }(i,j)\in B
\ee
\be
\al_i +\beta_j \leq  c_{ij}\;\mbox{ if }(i,j)\not\in B
\ee

\noindent Now we prove the following theorem.
\bt
Consider the problem P. Let $x_1\leq x_2\leq \cdots \leq x_n$ and $y_1\leq y_2\leq \cdots \leq y_n$
be given numbers and $c_{ij}=f(x_i-x_j) $ where $f:\;R\raro R$ is  convex. Then the North West corner solution is optimal for Problem P.
\et
\NI{\bf Proof. }
\NI Suppose $i < r$ and $j < s.$
Let $x_1\leq x_2\leq \cdots\leq x_i \leq \cdots\leq x_r\leq \cdots\leq x_n$ and $y_1\leq y_2\leq \cdots \leq y_j \leq \cdots\leq y_s\leq \cdots\leq y_n$ be given numbers.
It is shown that $$x_i-y_s \leq x_i-y_j\leq x_r-y_j$$
$$x_i-y_s \leq x_r-y_s\leq x_r-y_j$$
\noindent Then there exist $0\leq \lmd\leq 1$ and $0\leq \mu\leq 1$ such that
$$x_i-y_j = \lmd (x_i-y_s) + (1-\lmd) (x_r-y_j)$$
$$x_r-y_s = \mu (x_i-y_s) + (1-\mu)(x_r-y_j)$$
\begin{center}
	\begin{table}[htp]
		\caption{\scriptsize North West corner rule solution}
		{\scriptsize
			\begin{center}
				\begin{tabular}{|c|c|c|c|c|c|c|c|c|c|}
					\hline
					& $\cdots$&$j$ & $\cdots$ & $k$& $\cdots $& $l$ & $\cdots$    & $s$ & $\cdots$ \\ \hline
					$\vdots$ &&  & & &  &  &  & &\\\hline
					$i$ & & $*$ & & &  &  & & $(i,s)$ & \\\hline
					$\vdots$ & & $\vdots$ & & & & &  & & \\\hline
					$p$& &$*$ & $\cdots$ &$*$ & & &  & & \\\hline
					$\vdots$& & &  &$\vdots$ & &&  & &\\\hline
					$q$& &&  &$*$ &$\cdots$ & $*$&  & &\\\hline
					$\vdots$& &&  && & $\vdots$&  & & \\\hline
					$r$& &&  && & $*$&$\cdots$  &$*$ & \\\hline
					$\vdots$ &&  & & &  &  &  & &\\\hline
				\end{tabular}
			\end{center}
		}
	\end{table}
\end{center}

However since $x_i-y_j + x_r-y_s$  $= (\lmd + \mu) (x_i-y_s) + (2-\lmd-\mu)(x_r-y_j),$ it follows that $\lmd +\mu=1.$
Therefore, by convexity of $f$
$$f(x_i-y_j)\leq \lmd f(x_i-y_s) + (1-\lmd) f(x_r-y_j)$$
$$f(x_r-y_s)\leq \mu f(x_i-y_s) + (1-\mu)f(x_r-y_j)$$
$$f(x_i-y_j) + f(x_r-y_s)\leq (\lmd + \mu) f(x_i-y_s)$$
$$ + (2-\lmd-\mu)f(x_r-y_j)$$
Using $\lmd +\mu=1$ and $c_{ij}=f(x_i-x_j)$ it follows that
$c_{ij} + c_{rs} \leq c_{is} +c_{rj} \mbox{ for all } i,j, r, s \mbox{ such that }$ $i < r \mbox{ and }j < s.$
By Lemma 2.1, it follows that the North West corner rule produces an optimal solution. \qed

\bcl
In problem P, suppose  $c_{ij}=(x_i-x_j)^{2}$ or $c_{ij}=|i-j|.$  Then  North West corner solution is an optimal solution for Problem P.
\ecl

\section{Transportation Problem and a Weighted Version of K$\ddot{\mbox{o}}$nig-Egerv$\acute{\mbox{a}}$ry Theorem}
We now consider the cardinality of a maximum matching and the size of a minimum vertex cover in a bipartite graph. Consider the entries of a matrix $A=(a_{ij})\in R^{n\times n}$ as {\it
	points} and a row or a column as a {\it line.} A set of points is said to be
{\it independent} if none of the lines of the matrix contains more than one
point in the set. Suppose $T$ is an independent set of points. Then an element of $T$ is
said to be an independent point.
K$\ddot{\mbox{o}}$nig-Egerv$\acute{\mbox{a}}$ry theorem (Egerv$\acute{\mbox{a}}$ry \cite{7},
K$\ddot{\mbox{o}}$nig \cite{8})  is stated as follows:
\bt
Let $S$ be a nonempty subset of points of a matrix $A=(a_{ij})\in
R^{n\times n},$ then the maximum number of independent points that can be
selected in $S$ is equal to the minimum number of  lines covering all
points in $S.$
\et
K$\ddot{\mbox{o}}$nig-Egerv$\acute{\mbox{a}}$ry theorem is used to obtain Hungarian algorithm and it is used   to prove
the finite convergence of the Hungarian method for linear
assignment  problem.

Suppose that $\bar{C}=[\bar{c}_{ij}]\in R^{n\times n}$ is a cost matrix of the
assignment  problem. We obtain a reduced cost matrix $\bar{C}^{'}$  of order $n$ by subtracting the smallest element in each row and then subtracting the smallest element in each column. Note that all the elements of $\bar{C}^{'}$ are non-negative and there is  at least one zero in every row and every column.
Recall that any two zero is said to be {\it independent} if they do not  lie in the same line.  Let $t$ be the number of independent zeros in the  reduced cost matrix $\bar{C}^{'}$ and $t\leq n.$ The K$\ddot{\mbox{o}}$nig-Egerv$\acute{\mbox{a}}$ry theorem states that maximum cardinality of an independent set of $0$'s is equal to minimum number of lines to cover all $0$'s.

In this section, we describe a weighted version of
K$\ddot{\mbox{o}}$nig-Egerv$\acute{\mbox{a}}$ry theorem and use it to provide a weighted version of  Hungarian method for solving a transportation problem.  This states minium cut - maximum flow theorem of Ford and Fulkerson in a different way. Accordingly, a standard transportation problem can also be written as a linear assignment problem as follows:

Let $\cs_{1}\,=\,\{1,2,\cdots,a_{1}\},$
$\ct_{1}\,=\,\{1,2,\cdots,b_{1}\},$
$$\dsp { \cs_{r}=\,\{\sum_{j=1}^{r-1}a_{j}\,+ 1,\sum_{j=1}^{r-1} a_{j}\,+
	2,\cdots,\sum_{j=1}^{r}a_{j}\}},\, 2\leq r\leq m \mbox{ and}$$
$$\dsp { \ct_{s}=\,\{\sum_{j=1}^{s-1}b_{j}\,+ 1,\sum_{j=1}^{s-1} b_{j}\,+
	2,\cdots,\sum_{j=1}^{s}b_{j}\}}, 2\leq s\leq n. $$
Let $\dsp {\eta=\sum_{i=1}^{m}a_{i}=\sum_{j=1}^{n}b_{j}.}$
Consider a linear assignment problem in which total number of machines is equal to
total number of jobs ($\eta$). Let  $C=(c_{ij})\in R^{m\times n}$ be the
cost matrix  of the transportation problem. We construct a cost matrix
$\tilde{C}=(\tilde{c}_{ij})$ by copying $C_{i\cdot},$ $a_{i}$ times for
$i=1,\cdots,m$ and $C_{\cdot j},$ $b_{j}$ times for $j=1,\cdots,n.$ Thus
$\tilde{C}_{p\cdot}=C_{l\cdot}\;\fal\;p\in \cs_{l},$
$\tilde{C}_{\cdot p}=C_{\cdot l}\;\fal\;p\in \ct_{l}.$ The
matrix constructed in this manner leads to a square cost matrix $\tilde{C}$
of order $\eta\times \eta$ for the linear assignment problem. So we arrive
at an equivalent assignment problem of the transportation problem. Now looking at the equivalent linear assignment problem, we observe to see that there are $mn$ blocks  in $\tilde{C}$ where $(ij)^{th}$ block is
of size  $a_{i}\times b_{j}$ consisting of identical entries $c_{ij}.$

We explore the possibility of extending the Hungarian method for
transportation problem using the original cost matrix $C$ of order $m\times n.$
Note that in $\tilde{C},$ $(ij)^{th}$ block of size $a_{i}\times b_{j}$
consisting of identical  entries $c_{ij}$  can be treated as a single entry  in
$C$ in the $(ij)^{th}$ position.  We provide a weight $a_{i}$ for the
$i^{th}$ row in $C$ and a weight $b_{j}$ for the $j^{th}$ column in $C.$
We now state  a weighted version of K$\ddot{\mbox{o}}$nig-Egerv$\acute{\mbox{a}}$ry theorem.

In this theorem we use the following terminology. The entries of a
matrix $C=(c_{ij})\in R^{m\times n}$ are called {\it blocks}. The $i^{th}$
row $C_{i\cdot}$  is known as a horizontal  line with weight $a_{i}$ and
$C_{\cdot j,}$  the $j^{th}$ column is a vertical line with weight $b_{j}.$
A set of blocks is said to be
{\it independent} if none of the lines of the matrix contains more than one
block in the set. Suppose $\Lambda$ is an independent set of blocks.
Then an element of $\Lambda$ is said to be an independent block.

\noindent Now we prove the following theorem.
\bt
If $\Omega$ is a nonempty subset of the blocks of a matrix $C,$
then the maximum number of independent blocks that can be selected in
$\Omega$ is equal to the lines with minimum total weight covering all the
blocks in $\Omega.$
\et
\NI{\bf Proof. } Note that in $C,$ a horizontal line with weight $a_{i}$
is equivalent to $a_{i}$ rows and  a vertical line with weight $b_{j}$ is
equivalent to $b_{j}$ columns in $\tilde{C}.$ Let $\tilde{\Omega}$  be a
nonempty subset of points of a matrix  $\tilde{C}.$ Now by Theorem 3.1,
the maximum number of independent points
that can be selected in $\tilde{\Omega}$ is equal to the minimum number of
lines covering all elements in $\tilde{\Omega}.$ Now in $\tilde{C},$
drawing $a_{i}$ horizontal lines is equivalent to drawing a horizontal line
with weight $a_{i}$ in $C.$ Similarly in $\tilde{C},$ drawing $b_{j}$
vertical lines  is equivalent to drawing a vertical line with weight
$b_{j}$ in $C.$ Since we do not distinguish between horizontal and
vertical lines, the maximum number of independent blocks that can be
selected in $\Omega$ is equal to the lines with minimum total weight
covering all blocks of $\Omega.$ \qed

The Hungarian method for the linear assignment problem  was developed by
Kuhn \cite{9} which has computational complexity $O(n^{4}).$ Lawler \cite{10} developed an
order $O(n^{3})$ version of the algorithm. Cechl$\acute{\mbox{a}}$rov$\acute{\mbox{a}}$ \cite{5} observed that in practical situations, it may be useful  to get an overall picture about all the optima as well and obtains a generalization of the Berge's
theorem.

We now apply weighted  K$\ddot{\mbox{o}}$nig-Egerv$\acute{\mbox{a}}$ry
theorem to get an weighted version of Hungarian method for solving transportation problems. Note that by weighted  K$\ddot{\mbox{o}}$nig-Egerv$\acute{\mbox{a}}$ry
theorem,  maximum number of independent zero blocks that can be selected
is equal to the lines with minimum total weight covering all the
blocks.  We describe the weighted version of Hungarian method  based on
weighted K$\ddot{\mbox{o}}$nig-Egerv$\acute{\mbox{a}}$ry theorem for
solving the transportation problem.

The basic steps of the weighted version of  Hungarian method
are same as Hungarian method. The termination rule is as in Theorem
3.2, i.e., the weights of the  lines drawn with minimum total weight is
equal to $\eta.$  Therefore, the proof of finite convergence also follows from Theorem 3.2. For the
sake of completeness we provide the basic steps of the weighted version of  Hungarian
method  for solving the transportation problem:
\begin{description}
	\item{Step 1:}  Get the  reduced cost matrix by  subtracting the smallest element in each row and then subtracting the smallest element in each column.
	
	\item{Step 2:}  Draw lines with minimum total weight to cover all zero blocks. Let
	the total weight of the lines drawn be $\zeta$.
	\item{Step 3:}
	If $\zeta =\eta$, optimal matrix has been reached. Get an optimal
	solution  by assigning flows through blocks having zero entries.
	If $\zeta <\eta$, find the minimum of the entries not covered by any
	line.  Let it be $\delta$. Subtract $\delta$ from all uncovered entries
	and  add $\delta$ to all entries covered by two lines. With the new matrix
	so obtained go to step 2.
\end{description}

\subsection{Illustrative Example}
We illustrate the basic steps with an example.

\bex {\em Consider the following transportation problem.
	\vsp

	\begin{center}
		{\scriptsize Table 1}
		\vsp
		
		{\scriptsize
			\begin{tabular}{|l|cccc|l|}\hline
				\multicolumn{1}{|l}{Origin} & \multicolumn{4}{|c|}{Destination} &
				\multicolumn{1}{c|}{$a_{i}$}\\
				& $D_{1}$ & $D_{2}$ & $D_{3}$ & $D_{4}$ & \\ \hline
				$O_{1}$ & 10 & 7 & 3 & 6 & 3\\
				$O_{2}$ & 1 & 6 & 8 & 3 & 5\\
				$O_{3}$ & 7 & 4 & 5 & 3 & 7\\ \hline
				$b_{j}$ & 3 & 2 & 6 & 4 &  \\ \hline
		\end{tabular} }
	\end{center}
	
	\NI{\bf (i) } Subtracting the smallest element in each row and then subtracting the smallest element in each column, we get the  reduced cost matrix as in Table 2.
	
	\begin{center}
		{\scriptsize Table 2}
		\vsp
		{\scriptsize
			\begin{tabular}{lllr}
				$7--$& $3--$ & $0--$&$3$\\
				$|$  & $|$    &      & $|$\\
				$0$&$4$ &$7$&$2$\\
				$|$  & $|$    &      & $|$\\
				$4$&$0$ &$2$&$0$
		\end{tabular} }
	\end{center}

	\NI{\bf (ii) } In Table 2, row 1, column 1, column 2 and column 4 are crossed out.
	The lines to cover all zeros with minimum weight (12) is shown in
	Table 2. The lines drawn with minimum total weight is not equal to $\eta=15.$
	The minimum of the not crossed out elements is subtracted from these
	elements and added to the elements which are on the intersection of two
	lines. Continuing in this manner we find the optimal assignment in another  2
	iterations.  See Table 3 and 4.
	
	\begin{center}
		{\scriptsize Table 3}
		\vsp
		{\scriptsize
			\begin{tabular}{lllr}
				$9--$& $5--$ & $0--$&$5$\\
				$|$  &     &      & \\
				$0$&$4$ &$5$&$2$\\
				$|$  &     &      & \\
				$4--$&$0--$ &$0--$&$0$
		\end{tabular} }
	\end{center}

	\begin{center}
		{\scriptsize Table 4}
		\vsp
		{\scriptsize
			\begin{tabular}{lllr}
				$11--$& $5--$ & $0--$&$5$\\
				&     &      & \\
				$0--$&$2--$ &$3--$&$0$\\
				&     &      & \\
				$6--$&$0--$ &$0--$&$0$
		\end{tabular} }
	\end{center}
	
	\noindent We obtain the optimal solution based on the following approach. Inspecting the rows of the final reduced cost matrix we see that the
	row 1 contains only one zero which occurs in the $3^{rd}$ column. Since
	$\min(a_{1}, b_{3})=a_{1}=3,$ we have $x_{13}=3.$ We cross out  row 1
	(since no assignment will be made further) and update $b_{3}$ by
	$b_{3}-a_{1}=3.$ Now inspecting the columns we see that entry in the
	column 1 and row  2  contains $0.$ Since
	$\min(a_{2}, b_{1})=b_{1}=3,$ we have $x_{21}=3.$ We cross out  column 1.
	Continuing in this manner we have $x_{32}=2,$ $x_{33}=3,$ $x_{24}=2$ and
	$x_{34}=2.$
}
\eex
\brem
Note that it is quite likely that there is no single zero in any row
and column. In this situation, we arbitrarily select rows and columns with
minimum number of zeros.
\erem
\section{Conclusion} This paper considers some structured transportation  problems which arise in  sample surveys and other areas of statistics. For these  transportation problems, optimal solution can be obtained by applying the North West corner rule. We look at a weighted version of K$\ddot{\mbox{o}}$nig-Egerv$\acute{\mbox{a}}$ry theorem and  the corresponding version of Hungarian method. This helps to find the optimal solution of the transportation problem similar to the way we find the optimal solution of the assignment problem.


\section*{Acknowledgement}
This work is carried out under the project on Optimization and
Reliability  Modelling of Indian Statistical Institute. The author
R. Jana is thankful to the Department of Science and Technology,
Govt. of India, INSPIRE Fellowship Scheme for financial support.
\bibliographystyle{plain}
\bibliography{bibfile}

\end{document}